\documentclass[11 pt]{article}

\usepackage{mathrsfs} 

\usepackage[all]{xy}

\newtheorem{tm}{Theorem}[subsection]
\newtheorem{lm}[tm]{Lemma}
\newtheorem{pr}[tm]{Proposition}
\newtheorem{rmk}[tm]{Remark}

\newtheorem{??}[tm]{Question}

\font\tenmsb=msbm10
\font\sevenmsb=msbm7
\font\fivemsb=msbm5

\newfam\msbfam
\textfont\msbfam=\tenmsb
\scriptfont\msbfam=\sevenmsb
\scriptscriptfont\msbfam=\fivemsb
\def\Bbb#1{{\fam\msbfam #1}}

\font\teneufm=eufm10
\font\seveneufm=eufm7
\font\fiveeufm=eufm5
\newfam\eufmfam
\textfont\eufmfam=\teneufm
\scriptfont\eufmfam=\seveneufm
\scriptscriptfont\eufmfam=\fiveeufm
\def\frak#1{{\fam\eufmfam\relax#1}}

\def\lorw{\longrightarrow}
\newcommand\n{\noindent}

\newcommand\ci{\cite}

\newcommand\rat{{\Bbb Q}}
\newcommand\comp{{\Bbb C}}

\newcommand\blacksquare{{\hspace*{\fill} $\fbox{}$}}

\newcommand{\phix}[2]{ \,^p\!{\cal H}^{#1}(#2)}

\newcommand{\im}{ \hbox{\rm Im} }

\newcommand{\ke}{ \hbox{\rm Ker} }

\newcommand{\be}{\begin{equation}}
\newcommand{\ee}{\end{equation}}

\title{Exchange between  perverse and weight filtration for the Hilbert schemes
of points of two surfaces}

\author{Mark Andrea A.  de Cataldo\thanks{
Partially supported by N.S.A. and N.S.F.}\,,  
Tam\'as Hausel\thanks{Supported by a Royal Society University Research Fellowship},
Luca Migliorini\thanks{Partially supported by PRIN 2007 project ``Spazi di moduli e teoria di Lie''}}
\date{}
\begin{document}
\maketitle

\begin{abstract}
We show that a natural isomorphism between the rational cohomology groups of the two zero-dimensional Hilbert schemes of $n$-points of  two surfaces,  the affine plane minus the
 axes and the cotangent bundle of an elliptic curve, exchanges
the weight filtration on the first set of  cohomology  groups  with the perverse Leray filtration 
associated with a natural fibration on the second set of cohomology groups. We discuss some associated hard Lefschetz
phenomena.
\end{abstract}

\tableofcontents

\section{Introduction}
\label{secintro}
\subsection{The main result}
\label{subsecmre}
The theory of mixed Hodge structures 
endows
the rational cohomology groups $H^*(Z,\rat)$ of a complex algebraic variety $Z$ with the  increasing weight filtration ${\mathscr W}_Z$.
On the other hand, given   
 a projective map  $f:Z \lorw Z'$ from a nonsingular algebraic variety, 
the theory of perverse sheaves (with middle perversity)
endows  the rational cohomology groups $H^*(Z,\rat)$ with the 
increasing    perverse Leray filtration ${\mathscr P}_Z$ (see \ci{decmigso3}, for example).
The filtration ${\mathscr P}_Z$  is analogous to the classical Leray filtration
${\mathscr L}_Z$
for the continuous map $f$. If $Z$ is nonsingular, then perverse Leray filtration
 is the counterpart in cohomology of a splitting,
in the
 derived category of sheaves on $Z'$,
of the direct image complex  $Rf_*\rat_Z$ into a direct sum of shifted intersection cohomology complexes supported on closed subvarieties of $Z'$. 
In general, the two filtrations differ and
one has a strict inclusion
${\mathscr L}_Z\subseteq {\mathscr P}_Z$. 
The perverse  Leray filtration admits a geometric description
  in terms of basic algebraic topology
constructions; see the main result in  \ci{decmigso3}. 

In this paper, we give an example of the following  remarkable phenomenon
 relating the weight and perverse filtrations for
 two distinct, yet related, varieties. 
Let $n\geq 0$ be an integer. Let $E$ be an elliptic curve and let $X:= T^{\vee}(E)\simeq E \times \comp$ the  total space of the cotangent bundle of $E$. The Hilbert scheme $X^{[n]}$ of $n$ points on $X$ is
 a $2n$-dimensional nonsingular variety admitting a proper flat map
 $h_n: X^{[n]} \to {\comp}^{(n)} \simeq \comp^{n}$
 of relative dimension $n$ onto the $n$-th symmetric product $\comp^{(n)}$ of $\comp$.
 Let $Y:= \comp^*\times \comp^*$ and $Y^{[n]}$ be the corresponding Hilbert scheme of $n$-points
 on $Y$. The weights of the  mixed Hodge structure on $H^*(Y^{[n]}, \rat)$ are even,
 so that it makes sense to define the halved weight filtration
 $_{\frac{1}{2}}{\mathscr W}_{Y^{[n]}, k}:={\mathscr W}_{Y^{[n]}, 2k}$
 on $H^* (Y^{[n]}, \rat)$.

The main result of this paper is Theorem \ref{czziam} which establishes that there is a natural
isomorphism of graded vector spaces $\phi^{[n]}: H^*(X^{[n]}, \rat) \simeq
H^*(Y^{[n]}, \rat)$ that exchanges the perverse Leray filtration for the map $h_n$
with   the halved weight filtration defined above:
\[
\phi^{[n]} \left( {\mathscr P}_{X^{[n]}} \right) = {_{\frac{1}{2}}{\mathscr W}_{Y^{[n]}}}.\]
 In words, the isomorphism $\phi^{[n]}$ sends a class in
 ${\mathscr P}_{X^{[n]},p}$ to a class of type $(p,p)$ for ${\mathscr W}_{Y^{[n]}}$.
Theorem \ref{atq2}  relates the hard Lefshcetz theorem on the products of symmetric products
of the curve $E$ with the relative hard Lefschetz theorem for the map
$h_n$ and with a ``curious" hard Lefschetz theorem on the cohomology
of $Y^{[n]}$.

In the paper \ci{thmhcv},  we have proved that the same ``exchange of filtration" phenomenon appears in the  case of  the moduli space of  degree one and  rank two  Higgs bundles
on a curve of genus $g \geq 2$, endowed with the associated Hitchin map,
and of the associated twisted character variety. See \S \ref{conj} for a precise statement of this result.

The example dealt with in this paper presents a striking
difference with respect to  the one treated in \ci{thmhcv}. In the latter case, due to N\^go's
support theorem, most of
 the perverse sheaves showing up in the decomposition theorem are  supported on all of the target space  of the Hitchin map. On the other hand, in the  case treated here,
  every stratum in $\comp^{(n)}$ of the map $h_n: X^{[n]} \to \comp^{(n)}$ contributes  several
   perverse sheaves showing up in the decomposition theorem.
 
 At the moment, we cannot explain the exchange of filtration
phenomena described above, beyond the fact that we can observe them. 
In \S \ref{conj}  we discuss some properties shared by
the example considered in this paper and the one treated in \ci{thmhcv}, and we speculate 
on the  possibility of a  more general statement regarding the phenomenon of exchange of filtrations.

\subsection{Notation}\label{notation}

We work over the field of complex numbers $\comp$ and with singular cohomology with 
rational coefficients $\rat$. The results hold with no essential changes over 
any algebraically closed field and with $\rat_\ell$-adic cohomology. A variety
is a separated scheme of finite type over $\comp$.

We employ freely the language of derived categories,  perverse sheaves  and the decomposition 
theorem  as well as
 the language of Deligne's  mixed Hodge structures (MHS); the reader may consult
  \ci{bbd}, the survey \cite{bams} and  the textbooks  \ci{dimca,iv,k-s,petstee}. 
  For the convenience of the reader we summarize our notation and terminology below.

Given a variety $Z$, we work with the full subcategory $D_Z$ of the derived category of the
category of sheaves of rational vector spaces on $Z$ given by those bounded complexes $K$
on $Z$
whose cohomology sheaves ${\mathcal H}^i(K)$ on $Z$ are constructible; a sheaf on $Z$ is constructible
if there is a partition $Z = \coprod Z_a$ of $Y$ given by locally closed subvarieties such that
the restriction
$F_{|Z_a}$ is locally constant for every $a$. We denote the $i$-th perverse cohomology sheaf of a complex $K$ on $Z$ by
$\phix{i}{K}$; it is a perverse sheaf on $Z$.
Given a map $f: Z \to Z'$ of algebraic varieties, we denote the derived direct image
functor $Rf_*$ simply by $f_*$ and the $i$-th direct image  functor  by $R^if_*$. 

A filtration $F$ on a vector space is a finite increasing filtration $ \ldots \subseteq F_{i} V \subseteq
F_{i+1} V \subseteq \ldots $; finite means that  $F_{i} V = \{0\}$ for $i \ll 0$ and $F_i =V$ for $i \gg 0$.  A filtration $F$ on $V$  has type $[a,b]$
if $F_{a-1} V = \{0\}$  and $F_b V=V$.

Given a variety $Z$, the weight filtration on the cohomology groups $H^d(Z, \rat)$ is denoted by 
${\mathscr W}_Z$.  A map $f: Z \to Z'$ endows  the cohomology groups $H^d(Z,\rat)$
with two distinct filtrations, the Leray filtration ${\mathscr L}_Z$  and the perverse Leray filtration ${\mathscr P}_Z$.

In this paper, we are  concerned with the Hilbert schemes of $n$ points  $X^{[n]}$ and $Y^{[n]}$
associated with the two complex surfaces $X:=T^{\vee} E \simeq E \times \comp$, the total space of the cotangent bundle
of an elliptic curve $E$, and $Y:=\comp^* \times \comp^*$. We shall consider a
certain natural proper map $h_n: X^{[n]} \to \comp^{(n)}$.

\section{The Hilbert scheme of a surface and its cohomology groups}
\label{bho1}

\subsection{The decomposition theorem for the Hilbert-Chow map $\pi_n:S^{[n]} \to S^{(n)}$}
\label{dthc}
Let $S$ be a nonsingular connected complex analytic surface $S$ and $n\geq 0$
be a non-negative integer. We refer
the reader   to \ci{douady,chow,semires, gotts, nakajilectures}  for background
and references on Hilbert schemes of surfaces.   

We denote by $S^{(n)}:= S^n/{\frak S}_n$ the $n$-th symmetric product of 
$S$, i.e. the quotient of $S^n$ by the obvious action of the $n$-th
symmetric group.  
A partition of  $\nu=\{\nu_1, \ldots, \nu_l\}$ of $n$
is an unordered  collection of positive integers such that
  $\nu_1 + \ldots +\nu_l =n$; the integer $l=l(\nu)$ 
is called the  length of $\nu$. 
A point $x \in S^{(n)}$  gives rise to a partition
$\nu = \nu (x)$, for $x$ admits a unique
representation as a formal sum $\nu_1 s_1 + \ldots + \nu_l s_l$,
with $\nu_i$ positive integers adding up to $n$, and $s_i \in S$ distinct.
The subset $S^{(n)}_\nu  \subseteq S^{(n)}$ of points yielding the same partition
$\nu$ is a locally closed, irreducible, nonsingular subvariety of $S^{(n)}$
and we have that 
the symmetric product  $S^{(n)}$ is the disjoint union over the set of partitions on $n$
of these subvarieties:  $S^{(n)} = \coprod_\nu S^{(n)}_{\nu}$.
A partition $\nu$ gives rise to
a new variety $S^{(\nu)}$ as follows: represent the partition $\nu$
as a symbol $1^{a_1} 2^{a_2} \cdots n^{a_n}$, where $a_i$ is the number of times
$i$ appears in $\nu$; the $a_i \geq 0$, the length $l(\nu) = \sum a_i$
and $n= \sum_i i\,a_i$; finally, define $S^{(\nu)} := \prod_{i} S^{(a_i)}$
to be the indicated product of symmetric product of $S$.
If we define, ${\frak S}_\nu:= \prod {\frak S}_{a_i}$, then
$S^{(\nu)} = S^{l(\nu)} /{\frak S}_\nu$.
There is a natural finite map  $r^{(\nu)}:S^{(\nu)}\to S^{(n)}$
with image the 
closure $\overline{ S^{(n)}_\nu }$ and the resulting map
$S^{(\nu)} \to  \overline{ S^{(n)}_\nu }$ is the normalization
of the image.

The 
Hilbert scheme $S^{[n]}$ of zero-dimensional length $n$ subschemes
of $S$ is a connected complex manifold of dimension $2n$ and, if 
$S$ is algebraic, then so is $S^{[n]}$.  
There is the $n$-th  Hilbert-Chow map 
$\pi_n: S^{[n]} \to S^{(n)}$ sending a scheme to its support, counting multiplicities; this map is proper and it is a resolution of singularities of the symmetric product.

In view of  \ci{semires}, $\S$2.5, by using the correspondences in $S^{(\nu)} \times_{S^{(n)}} S^{[n]}$
inside $S^{(\nu)} \times S^{[n]}$
the decomposition theorem for the map $\pi_n$ yields a canonical isomorphism
in  the category ${D}_{S^{(n)}}$
 \be
\label{dths}
\xymatrix{ \gamma^{[n]}_S \, : \, = 
\sum_\nu \gamma_S^{(\nu)} \, : \, 
\bigoplus_{\nu}\, r^{(\nu)}_* \rat_{S^{(\nu)}} [2l(\nu)] \ar[r]^{\hskip 2.1cm\simeq} & \pi_* \rat_{S^{[n]}} [2n]. 
}
\ee

\subsection{The MHS on  $H^*(S^{[n]},\rat)$}
\label{oto}
If $S$ is algebraic, then,
by using the compatibility (see \ci{gysin})  with MHS  of the constructions leading to
the isomorphism  (\ref{dths}), we obtain a canonical  isomorphism of MHS 
(recall that a Tate twist in cohomology  $(-i)$ increases the weights by $2i$):
\be
\label{tat1}
\xymatrix{
\gamma_S^{[n]} =\sum \gamma^{(\nu)}_S :  \bigoplus_{\nu} 
\left(H^{* - 2 [n - l(\nu)]} \left(S^{(\nu)}, \rat \right) 
(l(\nu)-n) \right)
 \ar[r]^{\hskip 3.3cm \simeq} &
 H^*\left(S^{[n]}, \rat\right)}.
\ee
The fact that the  two sides of  (\ref{tat1}) are isomorphic has been first proved in
\ci{gottso} by using the theory of mixed Hodge modules.

Given a partition $\nu$ of $n$, consider the mixed Hodge substructure
\be\label{hnu}
H^*_{\nu} \left(S^{[n]}, \rat \right) := \im \, \gamma_S^{(\nu)} \simeq H^{ * -2 [n-l (\nu) ] }
\left(S^{(\nu)}, \rat \right) (l(\nu) - n)
\ee
so that the isomorphism of MHS (\ref{tat1}) now reads  as the internal direct sum decomposition
\be\label{222}
H^*\left( S^{[n]}, \rat \right) = \bigoplus_{\nu} H^*_{\nu}\left( S^{[n]}, \rat \right).\ee

\subsection{The map $\phi^{[n]}$ induced by a diffeomorphism $S_2 \simeq S_1$}
\label{ttvv}
The canonical isomorphism (\ref{tat1}) has the following  simple consequence.
Let $S_1$ and $S_2$ be two nonsingular  surfaces
and 
\be\label{0101}\phi:H^*(S_1, \rat) \simeq H^*(S_2,\rat)\ee
be an  isomorphism
of graded vector spaces. By taking tensor products and invariants,
 the map  $\phi$ induces, for  every partition $\nu$, 
an isomorphism of graded vector spaces
$\phi^{(\nu)}: H^*(S_1^{(\nu)}, \rat) \simeq 
H^*(S_2^{(\nu)}, \rat)$. 

By using the isomorphisms (\ref{dths}), we define
the map
\be\label{isodaus}
\phi^{[n]}: =    \left(\gamma_{S_2}^{[n]}\right) \circ    \left(\sum_{\nu} \phi^{(\nu)} \right) \circ \left(\gamma_{S_1}^{[n]}\right)^{-1}  :
H^*\left(S_1^{[n]}, \rat\right) \simeq H^*\left(S_2^{[n]}, \rat\right)
\ee
which is an isomorphisms of graded vector spaces.

If the surfaces $S_i$ are algebraic and $\phi$ is an isomorphism
of MHS, then so is (\ref{isodaus}). However, in this paper we  use
this set-up in the case: $S_1=E\times \comp$ ($E$ an elliptic curve)
and $S_2= \comp^*\times \comp^*$ and $\phi = \Phi^*$, where $\Phi: S_2\simeq S_1$ is a diffeomorphism.
In this case,  due to the incompatibility of the weights, $\phi$ and $\phi^{[n]}$ cannot be 
 isomorphisms of MHS.

It is likely that the results in \ci{voisin}  imply that if
we have a diffeomorphism $S_2 \simeq S_1$ of  nonsingular algebraic surfaces,
then there is a diffeomorphism $S_2^{[n]} \simeq S_1^{[n]}$. At present, we do not know this
and we do not need it here.

\section{The surfaces $X$ and $Y$ and the
 filtrations $_{\frac{1}{2}}{\mathscr W}_{Y^{[n]}}$ and  ${\mathscr P}_{Y^{[n]}}$}
 \label{consd}
For the remainder of the paper, we fix $n \geq 0$, an elliptic curve $E$ and
 we set 
\[Y:= \comp^* \times \comp^*, \qquad X:= T^{\vee} E \simeq E \times \comp, \]
i.e. $X$ is  the total space of the cotangent (canonical) bundle
of $E$. The isomorphism above is well-defined up to multiplication
by a non-zero scalar. 

The two surfaces $X$ and $Y$ are noncanonically diffeomorphic:
 choose $E$ to be $\comp / \Gamma$ where $\Gamma$ is the lattice of Gaussian integers;
 then use polar coordinates to identify $X$ and $Y$. 
Let $\Phi: Y \simeq X$ be any diffeomorphism and set
 $\phi := \Phi^* : H^*(X,\rat) \simeq H^*(Y,\rat).$
 We are in the situation  of $\S$\ref{ttvv}.(\ref{0101}) so that,  for every $n\geq 0$,
 we
 obtain the linear isomorphism (\ref{isodaus}) of graded  vector spaces
 \be\label{cvbt}
 \xymatrix{
 \phi^{[n]} : H^* \left( X^{[n]}, \rat \right) 
 \ar[rr]^{\simeq} && H^* \left( Y^{[n]}, \rat \right).
 }\ee
 
 As it was observed in $\S$\ref{ttvv}, for $n \geq 1$, the two sides are never isomorphic
 as MHS. In particular, (\ref{cvbt}) does not preserve the weight filtrations.
 
Let us remark that each $H^d(X^{[n]}, \rat)$ is a pure Hodge structure of
weight $d$.
 Since $H^*(X,\rat) \simeq H^*(E,\rat)$ is an isomorphism
of MHS, we have that the same is true for $H^*(X^{(\nu)}, \rat) \simeq
H^*(E^{(\nu)}, \rat)$ for every partition $\nu$ of $n$.
In view of the splitting of MHS (\ref{222}), we have
the following canonical isomorphism
of MHS
\[
H^*(X^{[n]},\rat) \stackrel{(\ref{222})}= \bigoplus_\nu H^*_{\nu} (X^{[n]}, \rat) 
\simeq \bigoplus _\nu H^{* -2[n-l(\nu)]} (E^{(\nu)}, \rat) ( l(\nu) -n).\]
Since each $H^{d}(E^{(\nu)}, \rat)$ is pure of weight $d$,
we conclude that
each $H^d(X^{[n]}, \rat)$ is pure of weight $d$ as well.
In particular, the weight filtration  ${\mathscr W}_{X^{[n]}}$
on $H^*(X^{[n]}, \rat)$ is simply the filtration by cohomological degree and this
should be contrasted with Proposition \ref{mhshsc}.

\subsection{The halved weight filtration $_{\frac{1}{2}}{\mathscr W}_{Y^{[n]}}$  on $H^* ( (\comp^* \times \comp^*)^{[n]}, \rat)$ }
\label{wfr}

In this section, we first compute the MHS on $H^*(Y^{[n]}, \rat)$
and determine the weight filtration
${\mathscr W}_{Y^{[n]}}$  on $H^*(Y^{[n]}, \rat)$. We then observe
that ${\mathscr W}_{Y^{[n]}}$ has no odd weights
so that we can define the halved weight  filtration $_{\frac{1}{2}}{\mathscr W}_{Y^{[n]},\,k}:= {\mathscr W}_{2k}$
on $H^*(Y^{[n]},\rat)$ by simply halfing the weights.

Recall that:  an MHS is of
Hodge-Tate type if the odd graded pieces of the weight
filtration are zero and  every even graded piece  ${\rm Gr}^{\mathscr W}_{2p}$
is of pure type $(p,p)$; an MHS is split of Hodge-Tate type if it is isomorphic
to a direct sum of pure MHS of Hodge-Tate type.

\begin{lm}
\label{mhscscs}
For every partition  $\nu $ of $n$,  the natural  MHS   on $H^*(Y^{(\nu)},  \rat)$ is split
of  Hodge-Tate type and, more precisely, 
\[
H^d\left(Y^{(\nu)}, \rat \right) \hbox{is  pure of weight $2d$ and Hodge-type  $(d,d)$},
\]
\[
0 = {\mathscr W}_{2d-1}  \subseteq {\mathscr W}_{2d} = H^d
 (Y^{(\nu)}, \rat).
\]
\end{lm}
{\em Proof.} Since $H^d(\comp^*,\rat)$ has type $(d,d)$,  for $d=0,1$,
and it is trivial otherwise, the statement follows 
from the K\"unneth isomorphism and the naturality of the mixed Hodge structure
for the inclusion $H^d(Y^{(\nu)},\rat)\subseteq H^d(Y^{l(\nu)},\rat)$ coming from the quotient map
$Y^{l(\nu)} \lorw Y^{l(\nu)}/{\frak S}_\nu = Y^{(\nu)}$.
\blacksquare

The following proposition is an immediate consequence of  Lemma \ref{mhscscs}
and of the equality of MHS (\ref{222}).

\begin{pr}
\label{mhshsc}
The natural mixed Hodge structure on $H^*(Y^{[n]},\rat)$
 is split of Hodge-Tate type. More precisely, in terms of the decomposition  
 {\rm (\ref{222})}, 
we have:
\[
{\mathscr W}_{2k}\,H^d \left(Y^{[n]},\rat\right)= 
\bigoplus_{d- (n-l(\nu))\leq k}  H^{d}_{\nu}
\left(Y^{[n]},\rat \right), \qquad {\mathscr W}_{2k}={\mathscr W}_{2k+1}.
\]
\end{pr}

Proposition \ref{mhshsc} allows us to define the halved weight 
filtration  $_{\frac{1}{2}}{\mathscr W}_{Y^{[n]}}$   by setting
\[
_{\frac{1}{2}}{\mathscr W}_{Y^{[n]},k}  := {\mathscr W}_{Y^{[n]}, 2k}.\]
The halved weight filtration $_{\frac{1}{2}}{\mathscr W}_{Y^{[n]}}$ on $H^*(Y^{[n]}, \rat)$  has type $[0,2n]$.

\subsection{Decomposition theorem for the Hitchin-like fibration
$h_n: X^{[n]} \to \comp^{(n)}$}
\label{anhif}

Let $p : X \to \comp$ be the induced projection.
Recall the notation in \S\ref{oto}. We have the commutative diagram
\be
\label{cdpe}
\xymatrix{
&     &   X^{[n]} \ar[d]^{\pi_n}   \ar@/^2.5pc/[dd]^{h_n}  \\
X^{l(\nu)}= E^{l(\nu)} \times \comp^{l(\nu)} \ar[r]^{\hskip 1.4cm / {\frak S}_{\nu}}   \ar[d]^{p^{l(\nu)}}   &
X^{(\nu)} \ar[r]^{r_X^{(\nu)}}   \ar[d]^{p^{(\nu)}}   &   X^{(n)} \ar[d]^{p^{(n)}} \\
\comp^{l(\nu)} \ar[r]^{/ {\frak S}_{\nu}} &\comp^{(\nu)} \ar[r]^{r_\comp^{(\nu)}}      &   \comp^{(n)}.
}
\ee
The maps $p^{l(\nu)}$  and $p^{(\nu)}$ are   of relative dimension $l(\nu)$ and the map $p^{(n)}$ is of relative dimension
$n$. In particular,
note that
\[
\dim  
\left\lbrace 
   { p^{ (n) } }^{-1} \left( \comp^{ (n) }_{\nu} \right)    
\right\rbrace 
= l(\nu) + n, \qquad
\dim  \left\lbrace{p^{(\nu)}}^{-1} \left(\comp^{(\nu)}_{\nu} \right)\right\rbrace = 2l(\nu).
\]
The fiber  of $p^{(\nu)}$ over  the general point of $\comp^{(\nu)}$ is isomorphic to
$E^{l(\nu)}$.  All the other fibers are isomorphic to  quotients of $E^{l(\nu)}$
under the action of suitable, not necessarily normal,   subgroups groups of
the finite group ${\frak S}_{\nu}$. The fibers
over the points in the small diagonal in $\comp^{(\nu)}$ are all  isomorphic
to $E^{(\nu)}= E^{l(\nu)}/{\frak S}_{\nu}$ so that, by the compatibility
with MHS of 
 Grothendieck's theorem on the rational cohomology of quotient varieties,
we have a canonical  isomorphism of MHS
\be\label{tgb}
H^*\left(E^{(\nu)},\rat\right) =  H^*\left(E^{l(\nu)},\rat\right)^{{\frak S}_{\nu}}.
\ee

The map $h_n: X^{[n]} \to \comp^{(n)}$ is projective of relative dimension 
$n = \frac{1}{2}\dim{X^{[n]}} = \dim{\comp^{(n)} }$, flat by \ci{matsu}, Corollary to Theorem 23.1,
and, as it is observed above, it has general fiber the Abelian variety $E^{l(\nu)}$.

\begin{rmk}\label{whywesay}{\rm
   We say that $h_n$ is a Hitchin-type map
because of the analogy it presents with the Hitchin map associated with the
moduli of Higgs bundles on a curve, where the dimensions of domain $M$, target $A$  and fibers $F$
are  related as above: $\dim{M}= 2\dim{A} = 2\dim{F}$ and also because
our main result  Theorem \ref{czziam} is analogous to the main result of
\ci{thmhcv}, which deals with rank two  Higgs bundles of odd degree on a curve.
}
\end{rmk}

Due to the commutativity of the diagram (\ref{cdpe}) and
to the functoriality of derived  push-forwards applied to $h_n = p^{(n)} \circ
\pi_n$, 
the decomposition theorem (\ref{dths}) for the map $\pi_n$  implies that
we have natural isomorphisms
\[
\xymatrix{
\bigoplus_{\nu}  { r^{(\nu)}_{\comp} }_{\hskip -0.2cm *}\,
p^{ (\nu) }_{*} \, \rat_{ X^{(\nu)} } [ 2l(\nu) ]
\ar[r]
& 
 \bigoplus_{\nu } 
p^{(n)}_* {r^{(\nu)}_X}_{\hskip -0.2cm *} \, \rat_{X^{(\nu)}} [2l(\nu)] \ar[r]^{
\hskip 0.8 cm \simeq} &
 {h_n}_* \rat_{X^{[n]}} [2n].
}
\]
By applying Grothendieck theorem on the invariant part 
of push-forwards under a quotient map under a finite  group action,
and by recalling that $p^{l(\nu)}$ is a projection map, we get a canonical
isomorphism
\[
p^{(\nu)}_* \rat_{X^{(\nu)}} = \left( p^{l(\nu)}_* \rat_{X^{(\nu)}}\right)^{{\frak S}_{\nu}}
= \bigoplus_{i=0}^{2l(\nu)} \left( R^i p^{l(\nu)}_* \rat \right)^ {{\frak S}_{\nu}} [-i].
\]
We thus get the distinguished splitting isomorphism in the category $D_{\comp^{(n)}}$
\be
\label{ftrg}
\xymatrix{
\Gamma^{[n]}_X: \bigoplus_{\nu}  
\bigoplus_{i=0}^{2l(\nu)}
\left\{  \left[{r^{(\nu)} _\comp}_{\hskip -0.2cm *} \,
  \left( R^i p^{l(\nu)}_* \rat \right)^ {{\frak S}_{\nu}}\right] \left[l(\nu)\right] \right\} \left[ -(i-l(\nu)) \right] \ar[rr]^{\hskip 4cm\simeq} && 
h_* \rat_{X^{[n]}} [2n].
 }
\ee
Since every  $r^{(\nu)}$ is finite, every  direct summand  in square brackets  is
an ordinary sheaf (not just  a complex).
Moreover,  since the  functors $r^{(\nu)}_*$ 
are $t$-exact,
 every  summand  in curly brackets is a perverse sheaf, in fact an intersection cohomology complex
 with twisted coefficients supported on  
$\overline{\comp^{(n)}_{\nu}}\subseteq \comp^{(n)}$.

It follows that (\ref{ftrg}) ``is" the decomposition theorem for the map $h_n$
in the sense that we decomposed the right-hand-side as  direct sum 
of shifted intersection cohomology complexes supported
on $\comp^{(n)}$. We note that, unlike the general statement of the decomposition 
theorem, we have obtained (\ref{ftrg}) as a distinguished isomorphism.

In order to simplify the notation, we set
\[
R^i_\nu : =
{r^{(\nu)} _\comp}_{\hskip -0.2cm *} \,
  \left( R^i p^{l(\nu)}_* \rat \right)^ {{\frak S}_{\nu}}.
\]

For our purposes, it is convenient to re-write (\ref{ftrg})
in the following two different ways, where the former emphasizes
the perverse-sheaf-nature of the summands, and the latter emphasizes the
ordinary-sheaf-nature of the summands. One merely needs to apply the appropriate shift
and re-organize the terms. By abuse of notation, we denote the resulting maps
with the same symbol $\Gamma^{[n]}_{X^{[n]}}$:
\begin{equation}\label{rvf}
\xymatrix{
\Gamma^{[n]}_X: 
 \bigoplus_{t=0}^{2n} \left(  \bigoplus_{i+ (n-l(\nu)) =t} R^i_\nu [l(\nu)] \right) [n-t]
 \ar[r]^{\hskip 2.3cm \simeq} & {h_n}_*\rat_{X^{[n]}} [n];
 }
 \end{equation}
 \begin{equation}\label{rvf2}
 \xymatrix{
\Gamma^{[n]}_X:  \bigoplus_{k=0}^{2n} \left( \bigoplus_{i+2(n-l(\nu)) =k} R^i_{\nu}   \right) [-k] \ar[r]^{\hskip 2.3 cm \simeq} &  {h_n}_*\rat_{X^{[n]}}.
 }
 \end{equation}

We now turn to the  decompositions in cohomology stemming from
the isomorphism(s)
$\Gamma^{[n]}_{X^{[n]}}$.
By taking components in (\ref{ftrg}), we have the equality of maps
in the derived category
 \[
 \Gamma^{[n]}_X= \sum_{\nu} \Gamma^{(\nu)}_X=
\sum_{\nu} \sum_{i=0}^{2 l(\nu)} \Gamma^{(\nu),i}_X\]
and, by taking the images in cohomology, we set
\be\label{gbbg}
{G}^{*}_{\nu}\left(X^{[n]},\rat\right) : = \im \, \Gamma_X^{(\nu)} \subseteq
H^* \left(X^{[n]}, \rat\right).
\ee

By  the very construction of the splitting (\ref{ftrg}), i.e. the
fact that is it obtained by pushing forward (\ref{dths}), we have that 
\[
G_{\nu} \left(X^{[n]}, \rat\right) :\stackrel{(\ref{gbbg})}=\im \, \Gamma^{(\nu)}_X = \im \, \gamma^{(\nu)}_X \stackrel{(\ref{hnu})}=: H^*_{\nu} \left(X^{[n]}, \rat \right)
\subseteq H^*\left( X^{[n]}, \rat \right),
\]
or, in words, the  two  distinguished splittings  of $H^*(X^{[n]},\rat)$  into $\nu$-components
arising from the decomposition theorem
for the Hilbert-Chow map $\pi_n$ and for the Hitchin-like map $h_n$
coincide.

For every fixed partition $\nu$ of $n$ and for every $0\leq i \leq 2l(\nu)$, 
we set 
\[
{H}^{*}_{\nu, i} \left(X^{[n]},\rat\right) : = \im \, \Gamma_X^{(\nu),i} \subseteq
H^*_{\nu} \left(X^{[n]}, \rat\right),
\]
so that
  \be\label{ddfr}
  H^*_{\nu} 
  \left(X^{[n]}, \rat \right) = \bigoplus_{i=0}^{2l(\nu)} H^*_{\nu, i} 
  \left(X^{[n]},\rat \right).\ee
  
  The following lemma shows that in each cohomological degree
 $d$, there is at most one non zero summand  $H^d_{\nu,i} (X^{[n]}, \rat)$
 in (\ref{ddfr}).
 
  \begin{lm}
  \label{saqew}
  We have the following
  \[
  H^q \left(\comp^{(n)}, R^i_{\nu}\right) \simeq \left\{ 
  \begin{array}{lr}
  0& \mbox{if}  \;\; q \neq 0, \\
  H^i \left(E^{(\nu)}, \rat \right)& \mbox{if} \;\; q=0.
  \end{array} \right.
  \]
  
 In particular, for every $d$, we have that
 \[
  H^d_{\nu} 
  \left(X^{[n]}, \rat \right) =  H^d_{\nu, i= d- 2(n -l(\nu))} 
  \left(X^{[n]},\rat \right) \simeq 
  H^{d- 2 (n -l(\nu)) } \left( E^{(\nu)}, \rat \right).
  \]
  \end{lm}
  {\em Proof.} We have
 \[
 {H}^{q} \left( \comp^{(n)},  
 {r^{(\nu)} _\comp}_{\hskip -0.2cm *} \,
  {\left( R^i p^{l(\nu)}_* \rat \right)}^{{\frak S}_{\nu}} \right) =
  H^{q} \left( \comp^{(\nu)}, \left(R^i p^{l(\nu)}_* \rat \right)^{{\frak S}_{\nu}} \right) 
  = H^{q} \left( \comp^{l(\nu)}, R^i p^{l(\nu)}_* \rat \right)^{{\frak S}_{\nu}} 
\]
Since $\comp^{l (\nu)}$ is contractible, the groups above 
  are zero whenever
  $q \neq 0$. In view of ({\ref{tgb}), for $q=0$ we have:
 \[
 H^0 \left( \comp^{l(\nu)}, \left(R^i p^{l(\nu)}_* \rat \right) \right)^{{\frak S}_{\nu}}
 = H^i \left(E^{l(\nu)}, \rat \right)^{{\frak S}_{\nu}}= 
 H^i \left(E^{(\nu)}, \rat \right).
 \]
 This proves the first statement.

  \n
According to (\ref{ftrg}) and the diagram (\ref{cdpe}),  
each summand  $H^d_{\nu, i} (X^{[n]},\rat)$ is the subspace of $H^d(X^{[n]}, \rat)$
 injective image of 
 \[
 \xymatrix{
 {H}^{d-2(n-l(\nu)) - i} \left( \comp^{(n)},  
 {r^{(\nu)} _\comp}_{\hskip -0.2cm *} \,
  {\left( R^i p^{l(\nu)}_* \rat \right)}^{{\frak S}_{\nu}} \right)
  & =&
  H^{d-2(n-l(\nu)) - i} \left( \comp^{(\nu)}, \left(R^i p^{l(\nu)}_* \rat \right)^{{\frak S}_{\nu}} \right) \\
  & =&\left(H^{d-2(n-l(\nu)) - i} \left( \comp^{l(\nu)}, \left(R^i p^{l(\nu)}_* \rat \right) \right)\right)^{{\frak S}_{\nu}} 
  } 
\]
The second statement now follows from (\ref{ddfr}) and from the first statement.
\blacksquare

 Summarizing: we have that for every $d$: 
 \be\label{bgr} 
  H^d\left(X^{[n]}, \rat\right) = \bigoplus_\nu H^*_{\nu} \left(X^{[n]}, \rat\right) =
  \bigoplus_\nu H^d_{\nu, d-2 (n - l(\nu))} \left( X^{[n]}, \rat\right) \simeq
  \bigoplus_\nu H^{d-2 (n - l(\nu)) } \left(E^{(\nu)}, \rat\right).
  \ee

\subsection{The   perverse Leray filtration
${\mathscr P}_{X^{[n]}}$ on $H^*(X^{[n]}, \rat)$}
\label{modplf}
The  theory of perverse sheaves endows $H^*(X^{[n]}, \rat)$
with the perverse Leray filtration  ${\mathscr P}_{X^{[n]}}$, i.e. with the perverse filtration
associated with the complex ${h_n}_* \rat_{X^{[n]}}[n]$; see \ci{decmigso3}.
Note that if we replace ${h_n}_* \rat_{X^{[n]}}[n]$  with
another shift ${h_n}_* \rat_{X^{[n]}}[m]$, the resulting filtrations gets translated.
We have made the choice  $m=n$ so that, in view of (\ref{rvf}),  the result has the 
same  type $[0,2n]$ as the one of $_{\frac{1}{2}}{\mathscr W}_{Y^{[n]}}$.

While, in general,  the perverse (Leray) filtration is canonically defined, there is no natural
splitting of it.
In our situation, in view of (\ref{rvf}) and of (\ref{bgr}), we have that
the perverse Leray filtration is naturally split:
\be\label{ggttrr}
{\mathscr P}_{X^{[n]},p} \, H^d \left( X^{[n]},\rat \right)= 
\bigoplus_{t \leq p} \bigoplus_{d - (n -l(\nu)) = t} H^d_{\nu} \left( X^{[n]},\rat \right)=  \bigoplus_{d- (n-l(\nu)) \leq  p } H^d_{\nu} \left( X^{[n]},\rat \right)   .
\ee

\begin{rmk}
\label{ler}
{\rm
In view of the expression (\ref{rvf2}),  it is straightforward to verify with the aid of
Lemma \ref{saqew} that the ordinary Leray filtration ${\mathscr L}_{X^{[n]}}$
on $H^*(X^{[n]}, \rat)$  for the map
$h_n$ is the filtration by cohomological degree. In particular, by comparing with (\ref{ggttrr}),
it is clear that the Leray filtration
is strictly included in the perverse Leray filtration. 
}
\end{rmk}

\section{The main result, relation with hard Lefschetz, and a speculation}
\label{mrconj}
  
 \subsection{``${\mathscr P}_{X^{[n]}} = {_{\frac{1}{2}}{\mathscr W}_{Y^{[n]}}}$"}
 \label{p=w}

 We are now ready to state and prove the main result of this paper.
 \begin{tm}
 \label{czziam}
 For every $n\geq 0$,
 the map $\phi^{[n]}$ {\rm (\ref{cvbt})} is a filtered isomorphism, i.e.
 \[
\phi^{[n]} ( {\mathscr P}_{X^{[n]}} ) = {_{\frac{1}{2}}{\mathscr W}_{Y^{[n]}}}.\]
 \end{tm}
 {\em Proof.} By its very definition, the map $\phi^{[n]}$ is a direct sum 
 map with respect to the $\nu$ decompositions   (\ref{222}) for  $S=X$ and $S=Y$, respectively
 It remains to apply Proposition \ref{mhshsc} and  (\ref{ggttrr}).
 \blacksquare

We would like to remark on the exceptional circumstance highlighted
by Theorem \ref{czziam}.
In view of the canonical splitting (\ref{ggttrr}),  
we say that a class $a \in H^d(X^{[n]}, \rat)$ has perversity $p$
if $a \in \oplus_{d-(n-l(\nu))=p}H^d(X^{[n]}, \rat)$.
Theorem \ref{czziam} shows that, regardless of the $(r,s)$ type
of $a$ with respect to  the pure Hodge structure 
$H^d(X^{[n]}, \rat)$, we have that $\phi^{[n]} (a) \in H^d(Y^{[n]}, \rat)$ has type
$(p,p)$ and, more precisely, lives in the $(p,p)$ part of the
split Hodge-Tate type structure.

The proof of Theorem \ref{czziam}  is heavily  based
on  the fact that we  have constructed  the  explicit   splitting (\ref{ggttrr}) of the perverse Leray filtration.
There is a different  approach  which is 
based on the following   geometric description  \ci{decmigso3} of  the perverse Leray
 filtration: 
let $s\geq 0$ and let $\Lambda^s \subseteq \comp^{(n)} \simeq \comp^n$ be 
a general  $s$-dimensional linear section of $\comp^n$;  then
$$
{\mathscr P}_{X^{[n]},p}\, H^d \left(X^{[n]},\rat \right)= \ke \,\left\{ 
H^d \left(X^{[n]},\rat \right) \lorw H^d \left(h_n^{-1}(\Lambda^{d-p-1}),\rat\right)\right\}.
$$ 
While we  omit the details of this approach, we do   point out the basic fact
leading to the  identification of  the kernel above with the right-hand-side of (\ref{ggttrr}):
 a general linear section $\Lambda^{d-p-1}$ avoids
the closure of a stratum $\overline{\comp_{\nu}^{(n)}}$, which has dimension
$l(\nu)$, if and only if $d - (n -l(\nu)) \leq p$.

\subsection{The curious hard Lefschetz (CHL) for $H^*( (\comp^* \times \comp^*)^{[n]}, \rat)$}
\label{chl}
Let $(z,w)$ be coordinates on $Y= \comp^*\times \comp^*$.
The $2$-form 
\[
\alpha_Y:= \frac{1}{(2 i \pi)^2} \frac{dz\wedge dw}{zw}
\]
is closed and defines an integral cohomology class which we denote with the same symbol. We have $\alpha_Y \in H^2(Y,\rat) \cap H^{2,2} (Y)$.
Let $p_i: Y^n \to Y$ be the $i$-th projection. 
Set $\alpha_{Y^n} = \sum_{i=1}^n p_i^* \alpha_Y \in H^2(Y^n,\rat) \cap H^{2,2} 
(Y^n)$.
Let $\alpha_{Y^{(n)}} \in H^2(Y^{(n)},\rat) \cap H^{2,2} (Y^{(n)})$
and $\alpha_{Y^{(\nu)}} \in H^2(Y^{(\nu)},\rat) \cap H^{2,2} (Y^{(\nu)})$
be the  naturally induced classes.
Let  $\alpha_{Y^{[n]}} := \pi_n^* \alpha_{Y^{(n)}}
\in H^2(Y^{(n)},\rat) \cap H^{2,2} (Y^{(n)})$
be the 
 pullback 
via the Hilbert-Chow map $\pi_n: Y^{[n]} \to Y^{(n)}$. 

Note that  because of Hodge type,  none of the $\alpha$-type classes
above is the first Chern 
 class of a holomorphic  line bundle on $Y^{[n]}$.
 Nonetheless,
  a simple explicit computation  based on Proposition \ref{mhshsc} shows that  cupping with the powers of  $\alpha_{Y^{[n]}}$,
gives rise to isomorphisms
\be
\label{chlanalogue}
\xymatrix{
{\rm Gr}_{2n-2k}^{{\mathscr W}_{Y^{[n]}}} H^{*}(Y^{[n]},\rat) \ar[rr]^{\alpha^k_{Y^{[n]}}}_\simeq&&
{\rm Gr}_{2n+2k}^{{\mathscr W}_{Y^{[n]}}} H^{* + 2k}(Y^{[n]},\rat).
}
\ee

These isomorphisms are  analogous to   the ``curious hard Lefschetz'' theorem of \ci{hauvill}.
Its curiosity consists of the fact that it is a statement concerning
a $(2,2)$ class on a noncompact variety, instead of  a $(1,1)$-class on a projective variety.
This apparently mysterious fact receives an explanation
from the coincidence 
of the halved weight filtration with the 
 perverse Leray filtration proved in the main Theorem \ref{czziam}.

\begin{??}\label{qn2}
{\rm
What corresponds to the CHL
(\ref{chlanalogue})
under the identification   $H^*(Y^{[n]}, \rat) \simeq H^* ( X^{[n]}, \rat)$
given by  (\ref{cvbt})? We answer this question in Theorem \ref{atq2}.
}
\end{??}

 \subsection{CHL on $Y^{[n]}$ $\Leftrightarrow$ the  HL on $E^{(\nu)}$ $\Leftrightarrow$ RHL for $h_n$}
 \label{chlrhl}
In this section,  we say that a rational cohomology class of degree two on a variety $Z$ 
is good (resp. ample)
if it is a non-zero   (resp. positive) rational multiple
of  the Chern class of an ample line bundle on $Z$.
The point of this definition is that the hard Lefschetz theorem
holds for a good class on a nonsingular projective manifold
as well as on its quotients by a finite group acting by algebraic
isomorphisms.

Fix any diffeomorphism $\Phi: Y=\comp^*\times \comp^* \simeq X=E\times \comp$. We obtain
the linear isomorphism (\ref{cvbt}) of graded vector spaces: $\phi_{[n]}: 
H^*(X^{[n]}, \rat) \simeq H^*(Y^{[n]}, \rat)$.

Let $\alpha_X, \alpha_{X^n}, \alpha_{X^{(n)}}, \alpha_{X^{(\nu)}}, \alpha_{X^{[n]}}$
be the classes obtained by transplanting the $\alpha$-classes defined starting from
$Y$ in section \ref{chl} via $\phi^{-1}_{[n]}$.

Note that by construction, for every surface $S$,
 the inclusion $H^*(S^{(n)}, \rat) \subseteq H^*(S^{[n]}, \rat)$ is given
 by the pull-back $\pi_n^*$ via the Hilbert-Chow map
 $\pi_n: S^{[n]} \to S^{(n)}$. In particular,
 we have that $\alpha_{X^{[n]}} = \pi^*_n \alpha_{X^{(n)}}$.
 This has to be verified in view of the fact
 that $\phi_{[n]}$ has not been defined using a diffeomorphism
 $Y^{[n]} \simeq X^{[n]}$
 between the Hilbert schemes.

Note that $\phi_{[n]}$ is not a map of MHS (this is already apparent for $n=1$).
 On the other hand, since $H^2(X,\comp) =H^2(E\times \comp)
\simeq H^2(E,\comp) = H^{1,1}(E)$, we see that
all the $\alpha$-classes $\alpha_X, \ldots, \alpha_{X^{[n]}}$ are in fact
in $H^2(-, \rat) \cap H^{1,1}(-)$. 
 
 Moreover,  the class $\alpha_X \in H^2(X,\rat) \simeq \rat$, being non-zero, is automatically good. In fact, it is ample if and only if the diffeomorphism
$\Phi: Y \simeq X$ preserves the  canonical orientations of the complex
analytic surfaces.

 It follows that
the $\alpha$-classes $\alpha_X, \alpha_{X^n}, \alpha_{X^{(n)}}$
and $\alpha_{X^{(\nu)}}$ are  good.  Since $\alpha_{X^{(\nu)}}$ is good, so is its restriction
to the fibers of $X^{(\nu)}\to \comp^{(\nu)}$.  The fibers of this map
 over points in the dense open stratum of $\comp^{(\nu)}$ consisting of multiplicity-free cycles
 are   isomorphic to the product  $E^{l(\nu)}$. Over the remaining points,
 the fibers are isomorphic to  finite quotients $E^{l(\nu)}/G$, where the   $G$ are suitable subgroups
of ${\frak S}_\nu$ (see section \ref{oto}).

On the other hand, if $n\geq 2$,
then $\alpha_{X^{[n]}}$ is not good: being a pull-back from $X^{(n)}$,
it is trivial on the positive dimensional projective fibers of the Hilbert Chow
birational map $\pi_n: X^{[n]}\to X^{(n)}$, a fact that prohibits 
goodness.

 In view of  the identifications of Lemma \ref{saqew} and of the fact that $\alpha_{X^{(\nu)}}$ and its restriction
 to $E^{(\nu)}$ are  good, 
 we have 
that the classical hard Lefschetz isomorphisms for the nonsingular
 projective  $E^{(\nu)}$ of dimension $l(\nu)$
 reads as follows
 \be
 \label{hlied}
  \alpha_{X^{(\nu)}}^j:
  H^{l(\nu) -j}_{\nu} \left(X^{[n]}, \rat \right)
 =
 H^{l(\nu) -j} \left(E^{(\nu)}, \rat
 \right)  
 \stackrel{\simeq}\lorw
 H^{l(\nu) +j} \left(E^{(\nu)}, \rat \right) = H^{l(\nu) +j}_{\nu} \left(X^{[n]}, \rat \right).
  \ee

  \begin{rmk}
  \label{actionalpha}
  {\rm
  Since $\alpha_{X^{[n]}}$ is a pull-back from $X^{(n)}$, 
  its action via cup product
  on ${\pi_n}_* \rat_{X^{[n]}} [2n]$ is diagonal with respect  to the decomposition
  into $\nu$-summands (\ref{dths}).
  Moreover,
  the induced action on each 
  $\nu$-summand
  is the action via cup product with $\alpha_{X^{(\nu)}}$. The same holds after taking cohomology.}
  \end{rmk}
  The hard Lefschetz isomorphisms (\ref{hlied}) express a property of this
  cup product action with  $\alpha_{X^{[n]}}$ in cohomology.
    In fact, (\ref{hlied})  is the reflection in cohomology
    of the fact that the conclusion
  of the  relative hard Lefschetz theorem (\ci{bbd}, Theorem 5.4.10; see also \ci{decmightam})
  holds 
  for the map $h_n: X^{[n]} \to \comp^{(n)}$ and  for the cup-product action with 
    $\alpha_{X^{[n]}}$, i.e. that 
  we have isomorphisms
  \be\label{rhlhnn}
  \alpha^j_{X^{[n]}} :  \phix{-j}{{h_n}_* \rat [2n]} \stackrel{\simeq}\lorw
  \phix{j}{{h_n}_* \rat [2n]},
  \ee
  where, in view of (\ref{ftrg}),  the perverse cohomology sheaves are
  \[
  \phix{j}{{h_n}_* \rat [2n]} =  \bigoplus_{i - l(\nu) =j} R^i_{\nu} [l(\nu)].\]
  In fact, the map of perverse sheaves (\ref{rhlhnn}) is defined
  simply because $\alpha_{X^{[n]}} \in H^2(X^{[n]}, \rat)$; see \ci{decmightam}), $\S$4.4.
  By using the identifications of Lemma
  \ref{saqew}, we deduce that the map (\ref{rhlhnn})
   is an isomorphism: in fact,  in view of  the isomorphisms (\ref{hlied}),  it is an isomorphism
  on the stalks of the respective cohomology sheaves.

  Recall that $\alpha_{X^{[n]}}$ is not good for $n\geq 2$, i.e. it is neither ``positive", nor "negative" on the fibers of $h_n$, so that the relative hard Lefschetz
  theorem does not apply in this context, yet we have (\ref{rhlhnn}).  This situation is similar to
  the one of the paper \ci{lef}, where the notion of lef line bundles has been introduced 
  and where it is proved that it is strongly linked to the hard Lefschetz theorem. The relation
  with the present situation is that, up to sign,  $\alpha_{X^{[n]}}$ is not ample on the fibers of $h_n$,
   but it is lef. 
  
  Recalling the expression (\ref{ggttrr}) for the perverse Leray filtration and
  Remark \ref{actionalpha}, a direct calculation using the hard Lefschetz isomorphisms
  (\ref{hlied}) and Theorem \ref{czziam}  implies  the following
  result, which answers Question \ref{qn2}.

  \begin{tm}
  \label{atq2}
  Under the identification $\phi_{[n]}: H^*(X^{[n]}, \rat) = H^*(Y^{[n]}, \rat)$, 
  the CHL {\rm (\ref{chlanalogue})}  becomes the (relative) hard Lefschetz
  {\rm (\ref{rhlhnn})}. 
  \end{tm}

We conclude this section by remarking that the splitting (\ref{ftrg}) of ${h_n}_* \rat_{X^{[n]}}$
has a remarkable property. Deligne's paper \ci{delignedercat} implies that
once we have the relative hard Lefschetz-type isomorphisms (\ref{rhlhnn}), we can construct three a priori 
distinct isomorphisms  between the l.h.s and the r.h.s of (\ref{ftrg}). 
Each one of these three splittings is characterized by a certain property of the matrices that express the
action 
of the cup product operations $\alpha_{X^{[n]}}^k: {h_n}_* \rat_{X^{[n]}} \to {h_n}_* \rat_{X^{[n]}} [2k]$ 
with respect to the splitting;
see \ci{delignedercat},
p.118 for the  definition of this matrix, Proposition 2.7 for the first splitting, section 3.1 for the second, and
Proposition 3.5 for the third. In general, these three splittings differ from each other,
e.g. in the case of  the projectivization of a  vector bundle with non trivial Chern classes, projecting over the base.

In our situation, there is the  fourth splitting (\ref{ftrg}).
The remarkable fact is that, in view of Remark \ref{actionalpha}, it is a matter of routine  to verify
that the  four splittings coincide.

\subsection{Speculating on where  to find the exchange of filtrations}
\label{conj}
The example treated in this paper and the one considered in \ci{thmhcv} have some properties in common which lead us to conjecture that the exchange of filtration occur for a certain
class of varieties and maps. Let us recall the main theorem of \ci{thmhcv}: 

Consider the  moduli space of semistable Higgs bundles ${\mathcal M}_{\rm Dol}$ parametrizing
stable rank 2 Higgs bundles  $(E,\phi)$ of degree $1$ on a fixed nonsingular projective curve
 $C$ of genus $g \geq 2$. There is  the Hitchin  proper and flat map $h: {\mathcal M}_{\rm Dol} \lorw \comp^{4g-3}$, which gives rise to the perverse Leray filtration
 ${\mathscr P}_{{\mathcal M}_{\rm Dol}}$.
By the {\em non-Abelian Hodge theorem},  ${\mathcal M}_{\rm Dol}$  is naturally diffeomorphic to the  twisted character variety 
$$
{\mathcal M}_{\rm B}:= \left\{ A_1,B_1,\dots, A_g,B_g \in {\rm GL}_2(\comp)\ | \
		 A_1^{-1} B_1^{-1} A_1 B_1 \dots A_g^{-1} B_g^{-1} A_g B_g = - {\rm I} \right\} / {\rm GL}_2(\comp)
$$ 
where the quotient is taken in the sense  of invariant theory.
The twisted character variety  ${\mathcal M}_{\rm B}$ carries a natural structure 
of nonsingular complex affine variety, with Hodge structure of Hodge-Tate type, with a natural splitting.

In complete analogy with Theorem \ref{czziam}, we have 
the main result   in \ci{thmhcv}, Theorem 4.2.9
\begin{tm}
In terms  of the isomorphism 
$H^*({\mathcal M}_{\rm B}) \stackrel{\simeq}{\lorw} H^*({\mathcal M}_{\rm Dol}) $
induced by the diffeomorphism ${\mathcal M}_{\rm B} \stackrel{\simeq}{\lorw} {\mathcal M}_{\rm Dol}$ stemming from the non-Abelian Hodge theorem, 
we have
$$
{\mathscr W}_{ {\mathcal M}_{B},2k } \, H^*({\mathcal M}_{\rm B})=
{\mathscr W}_{ {\mathcal M}_{B},2k+1 } \, H^*({\mathcal M}_{\rm B}) =
{\mathscr P}_{ {\mathcal M}_{\rm Dol},k} \, H^*({\mathcal M}_{\rm Dol}).
$$
\end{tm}

The varieties ${\mathcal M}_{\rm Dol}$ and $X^{[n]} $  
belong to the following class of varieties $Z$:

\begin{enumerate}
\item
$Z$ is a quasi-projective nonsingular variety  of even dimension $2m$ endowed with a  holomorphic 
symplectic structure $\omega\in H^0(Z;\Lambda^2T^*Z)$ and with  a $\comp^*$-action $\phi:\comp^*\times Z\to Z$, such that for $\phi_\lambda^*\omega=\lambda\omega$ for $\lambda\in \comp^*$ . 
\item
The ring $\Gamma(Z, {\mathcal O}_Z)$ is finitely generated and the affine reduction map  $h_Z: Z \lorw A={\rm Spec }\, \Gamma(Z, {\mathcal O}_Z)$
   is proper with  fibres of dimension $m$.
\item
The induced action on $A$ has a unique fixed point $o$ such that $\lim_{t \to 0}t \, y=o$ for all $y \in A$. 
\end{enumerate}

Let us note that, under these hypotheses, the Hodge structure on the cohomology groups
$H^d(Z,\rat)$ is pure of weight $d$:  the inclusion $h^{-1}(o) \subset Z$ induces an isomorphism
$H^d(Z, \rat) \simeq H^d(h^{-1}(o),\rat)$ of MHS; 
since  $Z$ is nonsingular,   the weight inequalities (\ci{hodge3} Theorem 8.2.4, iii. and iv.) imply the purity of $H^*(Z, \rat)$. 

Additionally we see that if $f$ and $g$ are functions in 
$\Gamma(A,{\mathcal O}_A)\cong \Gamma(Z,{\mathcal O}_Z)$ then we can write them as $f=\sum_{i>0} f_i$ and $g=\sum_{i>0} f_i$ and $g=\sum_{i>0} g_i$ such that $\phi_\lambda^*(f_i)=\lambda^if_i$ and $\phi_\lambda^*(g_i)=\lambda^ig_i$. Then the Poisson bracket satisfies 
$$\{f,g\}=\sum_{i,j>0}\{f_i,g_j\}=\sum_{i,j>0}\frac{1}{\lambda}\{\phi_\lambda^*f_i,\phi_\lambda^*g_j\}=\sum_{i,j>0}\lambda^{i+j-1}\{f_i,g_j\}.$$ Because $\lambda^k h =h$ for $k>0$ and generic $\lambda\in \comp^*$ only for the zero function, thus we can conclude $\{f,g\}=0$. 
Thus $h_Z$ is a completely integrable system. 

The two examples given in this paper and in \ci{thmhcv}   lead us to speculate
 whether it is  possible to associate with every variety $Z$ satisfying the 
 three assumptions above
another variety $\widetilde{Z}$ such that:

\begin{enumerate}
\item
$\widetilde{Z}$ is a quasi projective nonsingular variety endowed with a holomorphic 
symplectic structure.
\item
The affine reduction map  $h_{\widetilde{Z}}: \widetilde{Z} \lorw {\rm Spec }\, \Gamma(\widetilde{Z}, {\mathcal O}_{\widetilde{Z}})$   is birational (hence semismall in view of \ci{kale}, Lemma 2.11).
\item
There is a natural isomorphism $\phi: H^*(Z,\rat) \simeq H^*(\widetilde{Z}, \rat)$.
\item
The cohomology groups $H^*(\widetilde{Z},\rat)$ have a Hodge structure of split Hodge-Tate type.
\item 
Under the isomorphism $\phi$, the perverse filtration on $Z$ associated with the map $h$ corresponds to the halved weight filtration on $H^*(\widetilde{Z},\rat)$:
  a class of perversity $p$ on $Z$ would correspond
 to a class of type $(p,p)$ on $\widetilde{Z}$.

\end{enumerate}

Let us remark that, if the above were true, then the Hodge structure of $\widetilde{Z}$ cannot be pure.
In fact, in view of the relative hard Lefschetz
theorem, the class $\alpha \in H^2(Z, \rat)$ of any $h$-ample class  on $Z$ has necessarily perversity $2$. It would then follows that  $\phi (\alpha) \in  H^2(\widetilde{Z}, \rat)$ would have type $(2,2)$. 
In view of the conditions we have imposed on  the affine reduction maps of the two varieties,
i.e. the fact that $h_Z$ is a fibration with middle dimensional fibers and $h_{\widetilde{Z}}$
is semismall,
we like to think that $Z$ is ``as complete as possible," whereas $\widetilde{Z}$
is ``as affine as possible."  

At present, we do not know how to attack such a question and we still do not
know how to formulate  a principle that would justify the exchange of filtrations.

\end{document}